\documentclass[reqno, oneside]{amsart}
\usepackage{hyperref,url}
\usepackage{geometry}
\usepackage{pdfsync}
\usepackage{cite}

\theoremstyle{definition}

\newcommand\Atop[2]{\genfrac{}{}{0pt}{}{#1}{#2}}
\def\F(#1,#2,#3;#4){{}_2 F_1 \left ( \Atop{#1,#2} {#3} \, ; \, #4 \right )}
\def\Fthree(#1,#2,#3,#4,#5;#6){{}_3 F_2 \left ( \Atop{#1,#2,#3} {#4,#5} \, ; \, #6 \right )}
\def\f(#1,#2,#3;#4){F(#1,#2;#3;#4)}

\def\tfracx#1/#2@{\tfrac{#1}{#2}}

\usepackage[dvipsnames]{xcolor}

\begin{document}

\title{Elementary proof of congruences involving sum of binomial coefficients}
\markright{Binomialmod}
\author{Moa Apagodu}
\address{Virginia Commonwealth University\\ 1015 Floyd Ave, Richmond, 23284}
\maketitle

\begin{flushright}

To Olyad Apagodu 

\end{flushright}

\vskip .2in

\begin{abstract}
We provide elementary proof of several congruences involving single sum and multisums of binomial coefficients. 
\end{abstract}

\vskip .2in

\noindent {\bf INTRODUCTION}: We consider congruence of the form:
$$
\left(\sum_{k=0}^{rp^a-1} b(n+d)\right)\,\,\hbox{mod}\,\,\,p \,,
$$
where $b(n)$ is a combinatorial sequence, mainly binomial coefficients, $d\in \{0,1,2,3,4,\ldots,p^a\}$, $ a\in \{0,1,2,3,\ldots\}$, $r$ is a specific positive integer, and $p$ is an arbitrary prime. The case $a=1$ and single sum and single variable is covered in \cite{CHZ} and later extended to multisums and multivariables in \cite{AZ1}. This article generalizes the results in \cite{CHZ, AZ1} where the upper limit of summation is replaced by $rp^a-1$. 

\vskip .2in

\noindent The main ``trick'' in \cite{CHZ,AZ1} is {\bf The Freshman's Dream Identity} \cite{Wi}$(x+y)^{p} \equiv_p x^{p} + y^{p}$. In this article we use the same trick in the form $(x+y)^{p^a} \equiv_p x^{p^a} + y^{p^a}$. The proof follows from induction on $a$ and the case $a=1$. The second ingredient is {\bf Sum of a Geometric Series:} $\displaystyle{\sum_{i=0}^{n-1} z^i = \frac{z^n \, - \, 1}{z \, - \,1} }.$

\vskip .2in

\noindent {\bf NOTATION}: Let  $x \equiv_p y$ mean $x \, \equiv \, y \,\,(\hbox{mod} \,\, p)$, in other words, that $x-y$ is divisible by $p$.\\

\noindent The constant term of a Laurent polynomial $P(x_1,x_2,\ldots,x_n)$, alias
the coefficient of $x_1^0x_2^0..x_n^0$, is denoted by $CT[P(x_1,x_2,\ldots,x_n)]$.
The general coefficient of 
$x_1^{m_1}x_2^{m_2}..x_n^{m_n}$ in $P(x_1,x_2,\ldots,x_n)$ is denoted by 
$$
[x_1^{m_1}x_2^{m_2}..x_n^{m_n}]\,P(x_1,x_2,\ldots,x_n) \,.$$

\noindent {\bf Example 1.} 
$$
CT\left[ \frac{1}{xy}+3+5xy-x^3+6y^2\right]=3 \,, \,\,
[xy]\left[ \frac{1}{xy}+3+5xy+x^3+6y^2\right]=5.
$$

\noindent We use the symmetric representation of integers in $(-\frac{p}{2}, \frac{p}{2}]$ when reducing modulo a prime $p$. \\

\noindent {\bf Example 2.}   $6 \pmod{5}=1$ and $4 \pmod{5} = -1$. 

\vskip .2in 

\noindent We start by providing an elementary proof of, as advertised in the title, a result in \cite{ST}, where a complicated method is used to prove. 
 
\vskip .2in

\noindent {\bf Proposition $1$}[\cite{ST}, Remark 1.2] For any prime $p$ and $d\in \{0,1,2,3,\ldots,p^a\}$, we have
$$
\sum_{n=0}^{p^a-1} {2n \choose n+d} \,\,
\equiv_p  \,\,
\left\{\begin{array}{l}
1 ,\,\,\, \,\,\, \, \,\text{if} \, (\,p\,\equiv_3 \,2,  \,d\,\equiv_3 \,1 \,, a \,\,\text{is odd})\,\vee \,(\,p\,\equiv_3 \,1 \,\text{and}\,\,d\,\equiv_3 \,0) \vee\ (p\,\equiv_3 2\,,d\,\equiv_3 0\,, a \,\,\text{is odd})\\
-1 ,\,\,\,  \text{if} \, (\,p\,\equiv_3 \,2,  \,d\,\equiv_3 \,0 \,, a \,\,\text{is odd})\,\vee \,(\,p\,\equiv_3 \,1 \,\text{and}\,\,d\,\equiv_3 \,2) \vee\ (p\,\equiv_3 2\,,d\,\equiv_3 2\,, a \,\,\text{is even})\\
0,\,\,\,\,\,\,\,\,\,\text{otherwise}.
\end{array}\right.
$$

\begin{proof}

$$
\sum_{n=0}^{p^a-1} {2n \choose n+d} = \sum_{n=0}^{p^a-1} CT\left[ \left(\frac{(1+x)^{2n}}{x^{n+d}}\right)\right]
= \sum_{n=0}^{p^a-1} CT\left[ \left(2+x+\frac{1}{x}\right)^n \times \frac{1}{x^d}\right]
$$
$$
= CT\left[ \frac{\left(2+x+\frac{1}{x}\right)^{p^a}-1}{2+x+\frac{1}{x}-1} \times \frac{1}{x^d}\right]
\equiv_p 
CT \left[ \frac{2^{p^a}+x^{p^a}+\frac{1}{x^{p^a}}-1}{1+x+\frac{1}{x}} \times \frac{1}{x^d}\right] \, $$
$$
\equiv_p CT \left[ \frac{2+x^{p^a}+\frac{1}{x^{p^a}}-1}{1+x+\frac{1}{x}} \times \frac{1}{x^d}\right]
$$
$$
= CT \left[ \frac{1+x^{p^a}+\frac{1}{x^{p^a}}}{1+x+\frac{1}{x}}\times \frac{1}{x^d}\right]
\,= \,CT\left[ \frac{1+x^{p^a}+x^{2p^a}}{(1+x+x^2)x^{p^a-1}}\times \frac{1}{x^d}\right]
= [x^{p^a+d-1}]\, \left[ \frac{1+x^{p^a}}{1+x+x^2}\right]
$$
$$
= [x^{p^a+d-1}]\, \left[ \frac{(1-x)(1+x^{p^a})}{1-x^3}\right]\,
=[x^{p^a+d}]\, \left ( \, \sum_{i=0}^{\infty} x^{3i+1} \, \, + \, 
 \, \sum_{i=0}^{\infty} (-1) \cdot x^{3i+2} \, \right)+[x^{d}]\, \left ( \, \sum_{i=0}^{\infty} x^{3i+1} \, \, + \, 
 \, \sum_{i=0}^{\infty} (-1) \cdot x^{3i+2} \, \right).
$$
\noindent The result follows from extracting the coefficients of $x^{p^a+d}$ and $x^{d}$ in the above geometric series. For example, when $\,\,p\,\equiv \,2 (\hbox{mod} \,\, 3),\,\,  \,\,\,d\,\equiv \,2 (\hbox{mod} \, 3)$ and $a$ is even, then $\,\,\,p^a+d\,\equiv \,0 (\hbox{mod} \,\, 3) $ and the contribution from the first sum is zero. The only contribution is from the second sum, which is $-1$.
\end{proof}

\vskip .2in

\noindent Now we state and proof some of the results in \cite{CHZ} when the upper limit of summation is $rp^a-1$.

\vskip .2in

\noindent {\bf Proposition $1'$.} For any prime $p$ and $d\in \{0,1,2,3,4,5,\ldots,p^a\}$, we have
$$
\sum_{n=0}^{2p^a-1} {2n \choose n+d} \,\,
\equiv_p  \,\,
\left\{\begin{array}{l}
1 ,\,\,\, \,\,\, \, \,\text{if} \,\,p\,\equiv_3 \,1 \,\text{and}\,\,d\,\equiv_3 \,1\\
-4 ,\,\,\,   \text{if} \,\,p\,\equiv_3 \,1 \,\text{and}\,\,d\,\equiv_3 \,2\\

4 ,\,\,\,   \text{if} \,\,p\,\equiv_3 \,2 \,\text{and}\,\,d\,\equiv_3 \,1\\

-1 ,\,\,\,   \text{if} \,\,p\,\equiv_3 \,2 \,\text{and}\,\,d\,\equiv_3 \,2\\

3 ,\,\,\,   \text{if} \,\,p\,\equiv_3 \,1 \,\text{and}\,\,d\,\equiv_3 \,0\\

-3 ,\,\,\,   \text{if} \,\,p\,\equiv_3 \,2 \,\text{and}\,\,d\,\equiv_3 \,0\\
\end{array}\right.
$$

\begin {proof}
$$
\sum_{n=0}^{2p^a-1} {2n \choose n} = \sum_{n=0}^{2p^a-1} CT \left[ \left(  2+x+\frac{1}{x} \right)^n\times \frac{1}{x^d}\right]
=CT \left[ \frac{\left(2+x+\frac{1}{x}\right)^{2p^a}-1}{2+x+\frac{1}{x}-1}\times \frac{1}{x^d}\right]
$$
$$
=CT\left[ \frac{\left(6+4x+\frac{4}{x}+x^2+\frac{1}{x^2}\right)^{p^a}-1}{2+x+\frac{1}{x}-1} \times \frac{1}{x^d}\right]
\,\equiv_p \, CT\left[ \frac{\left(6+4x^{p^a}+\frac{4}{x^{p^a}}+x^{2{p^a}}+\frac{1}{x^{2{p^a}}}\right)-1}{2+x+\frac{1}{x}-1}\times \frac{1}{x^d}\right].
$$

\noindent Obviously only the terms $\displaystyle{5,\frac{4}{x^{p^a}}}$, and $\displaystyle{\frac{1}{x^{2p^a}}}$ contribute to the constant term. Discarding all the other ones, and simplifying, we get that this equals

\begin{eqnarray*}
[x^{2p^a+d-1}]\, \left[ \frac{1+4x^{p^a}+5x^{2p^a}}{1+x+x^2} \right] &=& [x^{2p^a+d-1}]\, \left[ \frac{1}{1+x+x^2}\right]  +4 \cdot [x^{p^a+d-1}]\, \left[ \frac{1}{1+x+x^2}\right] \\
& &+4 \cdot [x^{2p^a+d-1}]\, \left[ \frac{1}{1+x+x^2}\right]\\
&= & [x^{2p^a+d-1}]\, \left[ \frac{1-x}{1-x^3}  \right]  +4 \cdot [x^{p^a+d-1}]\, \left[ \frac{1-x}{1-x^3}\right] + [x^{d-1}]\, \left[ \frac{1-x}{1-x^3}\right]\\
&=&[x^{2p^a+d}]\, \left[ \sum_{i=0}^{\infty} x^{3i+1}  \right]  + [x^{2p^a+d}]\, \left[  \sum_{i=0}^{\infty} (-1) \cdot x^{3i+2} \right]  
+4 \cdot [x^{p^a+d}]\, \left[   \sum_{i=0}^{\infty} x^{3i+1} \right]  +\\
& & 4 \cdot [x^{p^a+d}]\, \left[   \sum_{i=0}^{\infty} (-1) \cdot x^{3i+2} \right]+5[x^{d}]\, \left[ \sum_{i=0}^{\infty} x^{3i+1}  \right]  + 5[x^{d}]\, \left[  \sum_{i=0}^{\infty} (-1) \cdot x^{3i+2} \right].
\end{eqnarray*}

\noindent The result follows from extracting the coefficients of $x^{p^a+d}$ and $x^d$ as in Proposition 2.
\end{proof}

\noindent The same method can be used to find the `` mod $p$" of $\displaystyle{\sum_{n=0}^{rp^a-1} { 2n \choose n+d}}$ for any specific $r$. 

\vskip .2in

\noindent Next we consider the Catalan numbers, $\displaystyle {C_n=\frac{1}{n+1}{ 2n \choose n}}$. 

\vskip .2in

\noindent {\bf Proposition $2.$} Let $C_n$ denote the $n$th Catalan number. Then, for any $p>3$, we have
$$
\sum_{n=0}^{p^a-1}  C_n \,\,
\equiv_p \,\,
\left\{\begin{array}{l}
1 ,\,\,\, \,\,\, \, \,\text{if} \,\,p\,\equiv 1 (\hbox{mod}\,\,  3) \,\,\text{or}\,\,p\,\equiv 2\, (\hbox{mod}\,\, 3)\,\,\text{and}\,\,a \,\,\text{is even}\\
-2 ,\,\,\, \text{if} \,\,p\,\equiv\, 2 (\hbox{mod}\,\, 3 )\,\,\text{and}\,\, a \,\,\text{is odd}.
\end{array}\right.
$$

\begin{proof}
Since  $C_n={{2n} \choose {n}} - {{2n} \choose {n-1}}$, it is readily seen that 
$$C_n=CT\left[(1-x)\left(2+x+\frac{1}{x}\right)^{n}\right].$$
We have
$$
\sum_{n=0}^{p^a-1}  C_n \,\,
= \sum_{n=0}^{p^a-1} CT \left[(1-x)\left(2+x+\frac{1}{x}\right)^n\right]
= CT \left[ \frac{(1-x)\left(\left(2+x+\frac{1}{x}\right)^{p^a}-1\right)}{2+x+\frac{1}{x}-1}\right]
$$
$$
\equiv_p CT\left[ \frac{(1-x)\left(\left(2+x^{p^a}+\frac{1}{x^{p^a}}\right)-1\right)}{2+x+\frac{1}{x}-1}\right]
(\hbox{By freshman's dream}).
$$
Since only the term $\frac{1}{x^{p^a}}$ in the numerator contributes to the constant term, this equals
$$
\, \, [x^{p^a-1}]\, \left[  \frac{1-x}{1+x+x^2}\right] \, = \,[x^{p^a-1}]\, \left[  \frac{(1-x)^2}{1-x^3}\right]
$$
$$
=  [x^{p^a}]\, \left[  \frac{x}{1-x^3}\right]
+[x^{p^a}]\, \left[  \frac{-2x^2}{1-x^3}\right]
+[x^{p^a}]\, \left[  \frac{x^3}{1-x^3}\right]
$$
$$
=  [x^{p^a}]\, \left[  \sum_{i=0}^{\infty} 1 \cdot x^{3i+1}  \right]
+[x^{p^a}]\, \left[    \sum_{i=0}^{\infty} (-2) \cdot x^{3i+2}   \right]
+[x^{p^a}]\, \left[   \sum_{i=0}^{\infty} 1 \cdot x^{3i+3} \right],
$$
and the result follows from extracting the coefficient of $x^{p^a}$ from the first or second geometric series above. (Note that we would never have to use the third geometric series, since $p>3.$)
\end{proof}

\noindent {\bf Proposition $2'.$} Let $C_n$ denote the $n$th Catalan number. Then, for every prime $p$,
$$
\sum_{n=0}^{2p^a-1}  C_n \,\,
\equiv_p \,\,
\left\{\begin{array}{l}
-7 ,\,\,\, \,\,\, \, \,\text{if} \,\,p\,\,\equiv \,2 (\hbox{mod}\,\,3)\,\,\text{and}\,\, a \,\,\text{is even}\\
\,\,\,\,2 ,\,\,\, \,\,\,\, \,\, \text{if} \,\,p\,\equiv \,\,1 (\hbox{mod}\,\,3)\,\,\text{or}\,\,p\,\equiv\, 2 (\hbox{mod}\,\, 3 )\,\,\text{and}\,\, a \,\,\text{is odd}.

\end{array}\right.
$$

\begin{proof}
Since  $C_n={{2n} \choose {n}} - {{2n} \choose {n-1}}$, it is readily seen that 
$$C_n=CT\left[(1-x)\left(2+x+\frac{1}{x}\right)^{n}\right].$$
We have
$$
\sum_{n=0}^{2p^a-1}  C_n \,\,
= \sum_{n=0}^{2p^a-1} CT \left[(1-x)\left(2+x+\frac{1}{x}\right)^{}\right]
= CT \left[ \frac{(1-x)\left(\left(2+x+\frac{1}{x}\right)^{2p^a}-1\right)}{2+x+\frac{1}{x}-1}\right]
$$

$$
= CT \left[ \frac{(1-x)\left(\left(6+4x+\frac{4}{x}+x^2+\frac{1}{x^2}\right)^{p^a}-1\right)}{2+x+\frac{1}{x}-1}\right]
$$

$$
\equiv_p CT \left[ \frac{(1-x)\left(\left(6+4x^{p^a}+\frac{4}{x^{p^a}}+x^{2p^a}+\frac{1}{x^{2p^a}}\right)-1\right)}{2+x+\frac{1}{x}-1}\right].
$$
Since only the terms $\frac{1}{x^{p^a}}$  and $\frac{4}{x^{p^a}}$ in the numerator contribute to the constant term, this equals
$$
\, \, [x^{2p^a-1}]\, \left[  \frac{(1-x)(1+4x^{p^a})}{1+x+x^2}\right] \, = \,[x^{2p^a-1}]\, \left[  \frac{(1-x)^2(1+4x^{p^a})}{1-x^3}\right]
$$
$$
=  [x^{2p^a}]\, \left[  \frac{x}{1-x^3}\right]
+[x^{2p^a}]\, \left[  \frac{-2x^2}{1-x^3}\right]
+[x^{2p^a}]\, \left[  \frac{x^3}{1-x^3}\right]
$$
$$
+  [x^{p^a}]\, \left[  \frac{4x}{1-x^3}\right]
+[x^{p^a}]\, \left[  \frac{-8x^2}{1-x^3}\right]
+[x^{p^a}]\, \left[  \frac{4x^3}{1-x^3}\right]
$$

$$
=  [x^{2p^a}]\, \left[  \sum_{i=0}^{\infty} 1 \cdot x^{3i+1}  \right]
+[x^{2p^a}]\, \left[    \sum_{i=0}^{\infty} (-2) \cdot x^{3i+2}   \right]
+[x^{2p^a}]\, \left[   \sum_{i=0}^{\infty} 1 \cdot x^{3i+3} \right],
$$

$$
+ [x^{p^a}]\, \left[  \sum_{i=0}^{\infty} 4 \cdot x^{3i+1}  \right]
+[x^{p^a}]\, \left[    \sum_{i=0}^{\infty} (-8) \cdot x^{3i+2}   \right]
+[x^{p^a}]\, \left[   \sum_{i=0}^{\infty} 4 \cdot x^{3i+3} \right],
$$

\noindent and the result follows from extracting the coefficients of $x^{2p^a}$ and $x^{p^a}$.
\end{proof}

\noindent The same method can be used to find the `` mod $p$" of $\displaystyle{\sum_{n=0}^{rp^a-1} C_n}$ for any specific $r$. 

\vskip .2in

\noindent The same method applied to the {\it Motzkin numbers}, $M_n$, that may be defined by the constant term formula
$$
M_n= CT \left [ \displaystyle{(1-x^2)\left(1+x+\frac{1}{x}\right)^n} \right ] ,
$$

\noindent leads to the following:\\

\noindent {\bf Proposition $3$.} Let $M_n$ denote the $n$th Motzkin number. Then, for any prime $p \geq 3$, we have
$$
\sum_{n=0}^{p^a-1}  M_n  \,
\equiv_p  \,\,
\left\{\begin{array}{l}
-2  ,\,\,\, \,\,\, \, \,\text{if} \,\,\,p\,\equiv \,1 (\hbox{mod}\,\,  4) \,\,\text{or}\,\,\,p\,\,\equiv \,\,3 (\hbox{mod}\,\, 4) \,\,\text{and}\,\,$a$\,\,\text{is even}\\
2 ,\,\,\,   \text{if} \,\,\,p\,\,\equiv \,\,3 (\hbox{mod}\,\, 4) \,\,\text{and}\,\,$a$\,\,\text{is odd}\end{array}\right.
$$

\begin {proof} 
$$
\sum_{n=0}^{p^a-1} M_n = \sum_{n=0}^{p^a-1} CT \left[(1-x^2)\left(1+x+\frac{1}{x}\right)^n\right]
= CT \left[ \frac{(1-x^2)\left(\left(1+x+\frac{1}{x}\right)^{p^a}-1\right)}{1+x+\frac{1}{x}-1}\right]
$$
$$
\, \equiv_p  \, CT \left[ \frac{(1-x^2)\left( 1+x^{p^a}+\frac{1}{x^{p^a}}-1 \right)}{1+x+\frac{1}{x}-1}\right]
\, = \,   CT \left[ \frac{(1-x^2)\left( x^{p^a}+\frac{1}{x^{p^a}} \right)}{x+\frac{1}{x}}\right]
\, = \,   CT \left[ \frac{x(1-x^2)\left( x^{p^a}+\frac{1}{x^{p^a}} \right)}{1+x^2}\right]
$$
$$
= [x^{p^a-1}]\, \left[  \frac{1-x^2}{1+x^2}\right]
= [x^{p^a}]\, \left[ \frac{x}{1+x^2}\right]  -  [x^{p^a}]\, \left[ \frac{x^3}{1+x^2}\right] 
$$
$$
= [x^{p^a}]\, \left[ \sum_{i=0}^{\infty} (-1)^i x^{2i+1} \right]  
+ [x^{p^a}]\, \left[   \sum_{i=0}^{\infty} (-1)^{i+1} x^{2i+3} \right] ,
$$
and the result follows from extracting the coefficient of $x^{p^a}$ from the first and second geometric series above, by noting that when $p \equiv \,1 \, (\hbox{mod}\,\,4)$ or $p \equiv \,3 \, (\hbox{mod} \,\,4)$ and $a$ is even, then $p^a\equiv \,1 \, (\hbox{mod}\,\,4) $ and in this case, $i$ is even in the first series, and odd in the second one, and vice-versa when  $p \equiv \,3 \, (\hbox{mod}\,\,4)$ and $a$ is odd.
\end{proof}

\noindent The same method can be used to find the `` mod $p$" of $\displaystyle{\sum_{n=0}^{rp^a-1} M_n}$ for any specific $r$. 

\vskip .2in

\noindent Next we consider the Apagodu-Zeilberger extension of Chen-Hou-Zeilberger method for discovery and proof of congruence theorems to multisums and multivariables\cite{AZ1} when the upper summation is replaced by $p^a-1$.

\vskip .2in

\noindent {\bf Proposition $4$.} For any prime $p$ and $a\in \{0,1,2,3,\ldots\}$, we have
$$
\sum_{n=0}^{p^a-1}\sum_{m=0}^{p^a-1} { n+m \choose m}^2\,\,
\equiv_p \,\,
\left\{\begin{array}{l}
1  ,\,\,\, \,\,\, \, \,\text{if} \,\,p\,\,\equiv \,\,1 (\hbox{mod}\,\, 3) \,\,\text{or}\,\,p\,\,\equiv \,\,2\, (\hbox{mod}\,\, 3)\,\,\text{and}\,\,\, a \,\,\text{is odd} \\
-1 ,\,\,\,\,p\,\,\equiv \,\,2\, (\hbox{mod}\,\, 3)\,\,\text{and}\,\,\, a \,\,\text{is even}\\
0\,\,\,\,\,\,\,\,\,p\,\,\,\,\equiv \,\,0\, (\hbox{mod}\,\, 3)\\
\end{array}\right.
$$

\begin{proof}
Let $\displaystyle{P(x,y)=\left(1+y\right)\left(1+\frac{1}{x}\right) \, }$ and  $\displaystyle{Q(x,y) = \left(1+x\right)\left(1+\frac{1}{y}\right) \, }$. Then 
$$
 { n+m \choose m}^2 =  { n+m \choose m} { n+m \choose n}=CT \left[ P(x,y)^n Q(x,y)^m\right].
$$
We have

\begin{eqnarray*}
\sum_{m=0}^{p^a-1}\sum_{n=0}^{p^a-1}  { m+n \choose m}^2  &=& \sum_{m=0}^{p^a-1}\sum_{n=0}^{p^a-1}  CT\left[ P(x,y) ^nQ(x,y)^m\right]\\
                                                                                             &=& CT\left[\sum_{m=0}^{p^a-1} \left[\frac{(P(x,y)^{p^a}-1)Q(x,y)^m}{P(x,y)-1}\right]\right]\\
                                                                                             &=& CT\left[\left(\frac{P(x,y)^{p^a}-1}{P(x,y)-1}\right) \left(\frac{Q(x,y)^{p^a}-1}{Q(x,y)-1}\right)\right] .
\end{eqnarray*}

\noindent We can pass to mod $p$ as above, and get

\begin{eqnarray*}
\sum_{m=0}^{p^a-1}\sum_{n=0}^{p^a-1}  { m+n \choose m}^2  &\equiv_p &CT\left[\left(\frac{P(x^{p^a},y^{p^a})-1}{P(x,y)-1}\right ) 
\left(\frac{Q(x^{p^a},y^{p^a})-1}{Q(x,y)-1}\right)\right]\\
&=& CT\left[\frac{(1+y^{p^a}+x^{p^a}y^{p^a})(1+x^{p^a}+x^{p^a}y^{p^a})}{(1+y+xy)(1+x+xy)x^{p^a-1}y^{p^a-1}}\right]\\
&=& [x^{p^a-1}y^{p^a-1}]\,\left[\frac{(1+y^{p^a}+x^{p^a}y^{p^a})(1+x^{p^a}+x^{p^a}y^{p^a})}{(1+y+xy)(1+x+xy)}\right]\\
&= & [x^{p^a-1}y^{p^a-1}]\,\left[\frac{1}{(1+y+xy)(1+x+xy)}\right] .
\end{eqnarray*}

\noindent It is possible to show that the coefficient of $x^ny^n$ in the Maclaurin expansion of the rational function $\frac{1}{(1+y+xy)(1+x+xy)}$ is $1$ when $n \equiv \,0 \pmod{3}$, $-1$ when $n \equiv \,1\,\pmod{3}$, and $0$ when $n \equiv \,2\,\pmod{3}$. 
One way is to do a partial fraction decomposition, and extract the coefficient of $x^n$, getting a certain expression in $y$ and $n$, and then extract the coefficient of $y^n$. Another way is by using the Apagodu$\--$Zeilberger algorithm (\cite{AZ2}), that outputs that the sequence of diagonal coefficients, let's call them $a(n)$, satisfy the recurrence equation $a(n+2)+a(n+1)+a(n)=0$,  with initial conditions $a(0)=1,a(1)=-1$. 
\end{proof}

\vskip .2in

\noindent We finally consider partial sums of {\it trinomial coefficients}.

\vskip .2in

\noindent {\bf Proposition $5$.}  Let $p>2$ be prime; then we have
$$
\sum_{m_1=0}^{p^a-1}\sum_{m_2=0}^{p^a-1} \sum_{m_3=0}^{p^a-1} { m_1+m_2+m_3 \choose m_1,m_2,m_3}  \, \equiv_p 1.
$$

\begin{proof}
First observe that 
${ m_1+m_2+m_3 \choose m_1,m_2,m_3}= CT \left[\frac{(x+y+z)^{m_1+m_2+m_3}}{x^{m_1}y^{m_2}z^{m_3}}\right]$.

Hence
$$
\sum_{m_1=0}^{p^a-1}\sum_{m_2=0}^{p^a-1} \sum_{m_3=0}^{p^a-1}{ m_1+m_2+m_3 \choose m_1,m_2,m_3}  = \sum_{m_1=0}^{p^a-1}\sum_{m_2=0}^{p^a-1} \sum_{m_3=0}^{p^a-1} 
CT \left[\frac{(x+y+z)^{m_1+m_2+m_3}}{x^{m_1}y^{m_2}z^{m_3}} \right]
$$
$$
= CT\left[\sum_{m_1=0}^{p^a-1}\sum_{m_2=0}^{p^a-1} \sum_{m_3=0}^{p^a-1}\frac{ (x+y+z)^{m_1+m_2+m_3}}{x^{m_1}y^{m_2}z^{m_3}}\right]
$$
$$
= CT \left[
\left ( \sum_{m_1=0}^{p^a-1} \left(\frac{x+y+z}{x}\right)^{m_1} \right )
\left ( \sum_{m_2=0}^{p^a-1} \left(\frac{x+y+z}{x}\right)^{m_2} \right )
\left ( \sum_{m_3=0}^{p^a-1} \left(\frac{x+y+z}{x}\right)^{m_3} \right )
\right] 
$$
$$
=CT\left[  \frac{(\frac{x+y+z}{x})^{p^a} -1}{\frac{x+y+z}{x}-1} \, \cdot \,
 \frac{(\frac{x+y+z}{y})^{p^a} -1}{\frac{x+y+z}{y}-1} \, \cdot \,
\frac{(\frac{x+y+z}{z})^{p^a} -1}{\frac{x+y+z}{z}-1} \right]
$$
$$
=[x^{p^a-1}y^{p^a-1}z^{p^a-1}]\, \left [ \frac{(x+y+z)^{p^a} -x^{p^a}}{y+z} \cdot \frac{(x+y+z)^{p^a} -y^{p^a}}{x+z}  \cdot
\frac{(x+y+z)^{p^a} -z^{p^a}}{x+y} \right ] .
$$
So far this is true for all $p$, not only $p$ prime. Now take it mod $p$ and get, using the freshman's dream in the form $(x+y+z)^{p^a} \equiv_p x^{p^a}+y^{p^a}+z^{p^a}$, that 
$$
\sum_{m_1=0}^{p^a-1}\sum_{m_2=0}^{p^a-1} \sum_{m_3=0}^{p^a-1} { m_1+m_2+m_3 \choose m_1,m_2,m_3}  
\,\, \equiv_p \,\,
[x^{p^a-1}y^{p^a-1}z^{p^a-1}]\, \left(\frac{y^{p^a}+z^{p^a}}{y+z} \cdot \frac{x^{p^a}+z^{p^a}}{x+z}\cdot \frac{x^{p^a}+y^{p^a}}{x+y}\right)
$$
$$
=[x^{p^a-1}y^{p^a-1}z^{p^a-1}]\, 
\left ( \sum_{i=0}^{p^a-1} (-1)^i y^i z^{p^a-1-i} \right)
\left ( \sum_{j=0}^{p^a-1} (-1)^j z^j x^{p^a-1-j} \right)
\left ( \sum_{k=0}^{p^a-1} (-1)^k x^{k} y^{p^a-1-k} \right)
$$
$$
=[x^{p^a-1}y^{p^a-1}z^{p^a-1}]\,  \left[
\sum_{0 \leq i,j,k<p^a} (-1)^{i+j+k} x^{p^a-1-j+k} y^{i+p^a-1-k} z^{p^a-1-i+j}
\right ] .
$$

\noindent The only contributions to the coefficient of $x^{p^a-1}y^{p^a-1}z^{p^a-1}$  in the above triple sum come when
$i=j=k$, so the desired coefficient of $x^{p^a-1}y^{p^a-1}z^{p^a-1}$
is
$$
\sum_{i=0}^{p^a-1} (-1)^{3i}= \sum_{i=0}^{p^a-1} (-1)^{i}=\underbrace{1-1+1-1+ \ldots +1-1}_{P^a-1}+ 1 = 1.
$$
\end{proof}

\noindent The same method of proof used in Proposition $5$ yields (with a little more effort) a {\it multinomial} generalization.\\

\noindent {\bf Proposition $7.$}   Let $p \geq 3 $ be prime,  then
$$
\sum_{m_1=0}^{p^a-1} \dots \sum_{m_n=0}^{p^a-1} 
{ m_1+\dots m_n \choose m_1, \dots , m_n}  \, \equiv_p 1.
$$

\vskip .2in

\noindent Finally, if we set $d=0$ in Corollary 1.1, equation 1.3, \cite{PS}, the right side simplifies to 

$$\sum_{n=0}^{p-1} (3n+1){ 2n \choose n}  \,\equiv_p\left\{\begin{array}{l}
-1 ,\,\,\, \,\,\, \, \,\text{if} \,\,\,p\,\equiv \,2 (\hbox{mod} \,\, 3)\\
1 ,\,\,\,   \text{if} \,\,p\,\equiv\,1 (\hbox{mod}\,\, 3)\,\,\,.\\
\end{array}\right.
$$
 
\vskip .2in

\noindent Motivated by this, we state two conjectures where the current method results in a rational function of higher degree that does not result the desired form.\\

\vskip .2in

\noindent {\bf Conjecture 1}: For any prime $p\geq 3$, we have 

$$\sum_{n=0}^{p-1}(5n+1) { 4n \choose 2n}  \,\equiv_p\left\{\begin{array}{l}
1 ,\,\,\, \,\,\, \, \,\text{if} \,\,\,p\,\equiv \,2 (\hbox{mod} \,\, 3)\\
-1 ,\,\,\, \,\,\, \, \,\,\,\text{if}\,\,\,p\,\equiv\, 1 (\hbox{mod}\,\, 3)\,\,\,.\\
\end{array}\right.
$$

\vskip .2in

\noindent The super Catalan numbers, first introduced by Ira Gessel \cite{IG}, also admits the following simple formulas. \\

\noindent {\bf Conjecture 2}: For any prime $p$, the supper Catalan number satisfies, 

$$\sum_{m=0}^{p-1}\sum_{n=0}^{p-1} \frac{{ 2m \choose m}{2n \choose n}}{{ n+m \choose n}}  \,\equiv_p\left\{\begin{array}{l}
1 ,\,\,\, \,\,\, \, \,\text{if} \,\,\,p\,\equiv \,1 (\hbox{mod} \,\, 3)\\
-1 ,\,\,\, \,\,\, \, \,\,\,\text{if}\,\,\,p\,\equiv\, 2 (\hbox{mod}\,\, 3)\,\,\,,\\
\end{array}\right.
$$

\vskip .2in

\noindent and 

$$\sum_{m=0}^{p-1}\sum_{n=0}^{p-1} (3m+3n+1) \frac{{ 2m \choose m}{2n \choose n}}{{ n+m \choose n}}  \,\equiv_p\left\{\begin{array}{l}
-7 ,\,\,\, \,\,\, \, \,\text{if} \,\,\,p\,\equiv \,\,1 (\hbox{mod}\,\,  3)\\
7 ,\,\,\, \,\,\, \, \,\,\,\text{if}\,\,p\,\equiv\, 2 (\hbox{mod}\,\, 3)\,\,\,.\\
\end{array}\right.
$$

\vskip .2in

\noindent {\bf Acknowledgements} We are grateful to Tewodros Amdeberhan and Doron Zeilberger for their valuable comments on an earlier version.

\end{document}